\documentclass[10pt]{article}
\usepackage{amssymb,amsmath,amsfonts,bbm,pifont,upgreek,bbold,accents}  %

\usepackage[colorlinks=true]{hyperref}
\hypersetup{urlcolor=blue, citecolor=red}

\font\teneufm=eufm10
\font\seveneufm=eufm7
\font\fiveeufm=eufm5
\newfam\eufmfam
\textfont\eufmfam=\teneufm
\scriptfont\eufmfam=\seveneufm
\scriptscriptfont\eufmfam=\fiveeufm

\newcommand\beq[1]{ \begin{equation}\label{#1} }
\newcommand{\eeq}{ \end{equation} }

\newcommand\beqa[1]{ \begin{eqnarray} \label{#1}}
\newcommand{\eeqa}{ \end{eqnarray} }
\newcommand{\beqano}{ \begin{eqnarray*} }
\newcommand{\eeqano}{ \end{eqnarray*} }

\newtheorem{theorem}{Theorem}
\newtheorem{definition}{Definition}
\newtheorem{proposition}{Proposition}
\newtheorem{lemma}{Lemma}
\newtheorem{sublemma}{Sublemma}
\newtheorem{remark}{Remark}
\newtheorem{notationalremark}{Notational Remark}
\newtheorem{corollary}{Corollary}
\newtheorem{assumption}{Assumption}
\newtheorem{claim}{Claim}

\newtheorem{tools}{$\negsp\negsp$}[subsection]

\newcommand\thm[1]{ \begin{theorem}\label{#1}}
\newcommand\thmtwo[2]{ \begin{theorem}[#1]\label{#2}}
\newcommand\ethm{ \end{theorem} }
\newcommand\dfn[1]{ \begin{definition}\label{#1} \rm}
\newcommand\dfntwo[2]{ \begin{definition}[#1]\label{#2} \rm}
\newcommand\edfn{ \end{definition} }
\newcommand\pro[1]{ \begin{proposition}\label{#1}}
\newcommand\protwo[2]{ \begin{proposition}[#1]\label{#2}}
\newcommand\epro{ \end{proposition} }
\newcommand\lem[1]{ \begin{lemma}\label{#1}}
\newcommand\lemtwo[2]{ \begin{lemma}[#1]\label{#2}}
\newcommand\elem{ \end{lemma} }
\newcommand\sublem[1]{ \begin{sublemma}\label{#1}}
\newcommand\sublemtwo[2]{ \begin{sublemma}[#1]\label{#2}}
\newcommand\esublem{ \end{sublemma} }
\newcommand\rem[1]{ \begin{remark}\label{#1} \rm}
\newcommand\erem{ \end{remark} }
\newcommand\notrem[1]{ \begin{notationalremark}\label{#1} \rm}
\newcommand\enotrem{ \end{notationalremark} }
\newcommand\cor[1]{ \begin{corollary}\label{#1}}
\newcommand\cortwo[2]{ \begin{corollary}[#1]\label{#2}}
\newcommand\ecor{ \end{corollary} }
\newcommand\asmp[1]{ \begin{assumption}\label{#1}}
\newcommand\asmptwo[2]{ \begin{assumption}[#1]\label{#2}}
\newcommand\easmp{ \end{assumption} }
\newcommand\clm[1]{ \begin{claim}\label{#1}}
\newcommand\eclm{ \end{claim} }
\newcommand\equ[1]{{\rm (\ref{#1})}}
\expandafter\chardef\csname pre amssym.def
at\endcsname=\the\catcode`\@
\catcode`\@=11
\def\undefine#1{\let#1\undefined}
\def\newsymbol#1#2#3#4#5{\let\next@\relax
 \ifnum#2=\@ne\let\next@\msafam@\else
 \ifnum#2=\tw@\let\next@\msbfam@\fi\fi
 \mathchardef#1="#3\next@#4#5}
\def\mathhexbox@#1#2#3{\relax
 \ifmmode\mathpalette{}{\m@th\mathchar"#1#2#3}%
 \else\leavevmode\hbox{$\m@th\mathchar"#1#2#3$}\fi}
\def\hexnumber@#1{\ifcase#1 0\or 1\or 2\or 3\or 4\or 5\or 6\or 7\or
8\or
 9\or A\or B\or C\or D\or E\or F\fi}
\ifcase\@ptsize
 \font\tenmsb=msbm10
 \font\sevenmsb=msbm7
 \font\fivemsb=msbm5
\or
 \font\tenmsb=msbm10 scaled \magstephalf
 \font\sevenmsb=msbm7 scaled \magstephalf
 \font\fivemsb=msbm5  scaled \magstephalf
\or
 \font\tenmsb=msbm10 scaled \magstep1
 \font\sevenmsb=msbm7 scaled \magstep1
 \font\fivemsb=msbm5 scaled \magstep1
\fi
\newfam\msbfam
\textfont\msbfam=\tenmsb
\scriptfont\msbfam=\sevenmsb
\scriptscriptfont\msbfam=\fivemsb
\edef\msbfam@{\hexnumber@\msbfam}
\def\Bbb#1{\fam\msbfam\relax#1}
\def\widehat#1{\setboxz@h{$\m@th#1$}%
 \ifdim\wdz@>\tw@ em\mathaccent"0\msbfam@5B{#1}%
 \else\mathaccent"0362{#1}\fi}
\def\widetilde#1{\setboxz@h{$\m@th#1$}%
 \ifdim\wdz@>\tw@ em\mathaccent"0\msbfam@5D{#1}%
 \else\mathaccent"0365{#1}\fi}

\def\RIfM@{\relax\ifmmode}
\def\nonmatherr@#1{\errmessage{\string#1\space allowed only in math mode}}
\def\Bbb{\RIfM@\expandafter\Bbb@\else
 \expandafter\nonmatherr@\expandafter\Bbb\fi}
\def\Bbb@#1{{\Bbb@@{#1}}}
\def\Bbb@@#1{\fam\msbfam\relax#1}
\def\setboxz@h{\setbox\z@\hbox}
\def\wdz@{\wd\z@}
\catcode`\@=\csname pre amssym.def at\endcsname
\newcommand{\giu}{{\medskip\noindent}}
\newcommand{\Giu}{{\bigskip\noindent}}
\newcommand{\nl}{{\smallskip\noindent}}

\newcommand{\qed}{\hskip.5truecm
\vrule width 1.7truemm height 3.5truemm depth 0.truemm
\par\Giu}

\newcommand{\negsp}{\hspace{-.09truecm}}  %%% equivalente a \!

\newcommand\ugper[1]{ \stackrel{#1}{=} }

\newcommand{\dst}{\displaystyle}

\newcommand\su[1]{ \frac{1}{ {#1}} }

\renewcommand{\natural}{ {\Bbb N}   }
\newcommand{\real}{ {\Bbb R}   }
\newcommand{\integer}{ {\Bbb Z}   }

\renewcommand{\a }{ {\alpha}   }

\newcommand{\g}{ {\gamma}   }

\newcommand{\x }{ {\xi}   }

\newcommand{\p}{ {\pi}   }

\renewcommand{\t}{ {\tau}   }

\renewcommand{\o}{ {\omega}   }

\newcommand\bks{\backslash}

\renewcommand\subset{\subseteq}

\title{\bf
Isolated Diophantine numbers
}

\begin{document}

\author{ 
F. Argentieri, L. Chierchia
%\\ \footnotesize Dipartimento di Matematica e Fisica
%\\ {\footnotesize @}
%\\ 
}

\maketitle

%{\small 
%\tableofcontents
%}
%%\date{}
%

\begin{abstract}

\noindent
In this short note, we discuss the topology of  Diophantine numbers, giving  simple explicit
examples of Diophantine isolated numbers (among those with same Diophantine constatnts), showing that, {\sl Diophantine sets are not always Cantor sets}. \\
General  properties of isolated Diophantine numbers are also briefly discussed.

\nl

\nl

\end{abstract}

\section{Introduction}
Diophantine numbers are irrational numbers badly approximated by rationals, namely, real numbers $\x$ satisfying, for some $\g,\t>0$,
\beq{D}
|\x q-p|\ge \frac{\g}{q^\t}\,,\qquad \forall \ p\in\integer\,, \ q\in \natural=\{1,2.,...\}\,.
\eeq
Such numbers, that form a set of full Lebesgue measure in\footnote{See, e.g., \cite{SM}, end of \S~25.} $\real$,  arise naturally in number theory and in small divisor problems in dynamics. Indeed, 
since the seminal works of C.L. Siegel, in the context of linearization of holomorphic diffeomorphisms around a fixed point \cite{S}, and of A.N.~Kolmogorov, in the context of Hamiltonian systems \cite{K}, Diophantine conditions as in \equ{D} (or higher dimensional analogs) are ubiquitous in perturbative Hamiltonian dynamics both in finite and infinite dimensions. 

\nl
Let us denote  by $D_{\g,\t}$ the Diophantine set of all real numbers satisfying condition \equ{D}
with fixed\footnote{Notice that $D_{\g,\t}=\emptyset$ whenever $\g\ge 1/2$ (trivially, taking $q=1$ and the minimum over $p$) or, (by Dirichlet's Theorem), when $\t<1$.} 
$0<\g<1/2$ and $\t\ge 1$. Clearly, $D_{\g,\t}$ is a closed and nowhere dense set. It is, therefore,  natural to ask whether $D_{\g,\t}$ is actually  a Cantor set, i.e., if it is also a perfect set (no isolated points). 

\nl
A (authoritative) place where the term Cantor set appears in association with Diophantine sets 
is Chapter III of the fundamental  and  beautiful book, {\sl Lectures on Celestial Mechanics},  by C.L. Siegel  and J.K.Moser \cite{SM}. In \S~32--36 of \cite{SM}, Moser, extending the previous text of Siegel\footnote{Compare with the 1971 Preface to the English Edition of \cite{SM}.}, includes a proof, in the analytic case, of his theorem on the persistence of invariant curves for area--preserving twist diffeomorphisms of the annulus \cite{M}. As well known, one of the main hypothesis is that $\o/2\p$ belongs to $D_{\g,\t}$ for given  $\g,\t$, where $\o$ denotes the rotation number of the unperturbed invariant curve. Moser call {\sl admissible} such numbers and, at page 245, writes\footnote{Obviously, Moser considers  Diophantine sets with exponent $\t>1$ (see p. 242), as it is well known that for $\t=1$, $D_{\g,\t}$ is not a Cantor set (see, e.g., \cite{B}).}:  ``{\sl the set of admissible values for $\o$ form a Cantor set of positive measure}''. Although, it does not appear a formal statement about the sets $D_{\g,\t}$, reading p. 245 of \cite{SM}, one might be lead to the belief that Diophantine sets are Cantor sets.

\nl
However, it turns out that, in general, this is {\sl not} the case: In \S~1, we show that the quadratic numbers $\a:=(n+\sqrt{n^2+4})/2=[n,n,n,...]=[\bar n]$ (in continued fraction expansion) are, for any $n\ge 2$, isolated in $D_{\g,\t}$ with 
$\g:=1/\a$ and $\t:={\log \a}/{\log n}$.

\nl
In \S~2, we  briefly review  some general properties of isolated Diophantine points, proven in \cite{F1,F2}, which show, in particular, that isolated Diophantine points are not that rare; \S~3 contains concluding remarks.

\section{Elementary examples of isolated Diophantine numbers}

{\bf Theorem 1}
{\sl
Let $n\in \natural$, $n\ge 2$ and define 
\beq{gt}
\a:=\frac{n+\sqrt{n^2+4}}2\ ,\ \qquad \g:=\frac1\a\ ,\qquad\ \t:=\frac{\log \a}{\log n}\ .
\eeq
Then, $\a$ is an isolated point of $D_{\g,\t}$.
}

\Giu
{\bf Remark}  (i) By  definition of $D_{\g,\t}$, it follows immediately that 
\beq{compl}
I_{\g,\t}(p,q):=\big\{\x\in\real:\, \big|\x -\frac{p}{q}\big|< \frac{\g}{q^{\t+1}}\big\}\subset \real\bks D_{\g,\t}\ ,\quad \forall q\in\natural, \forall p\in\integer\ .
\eeq
(ii) $D_{\g,\t}$ is invariant by translations by integers, as,
for $k\in\Bbb{Z}$,  $\xi\in D_{\g,\t}\iff \xi+k\in D_{\g,\t}$. Therefore, 
$\a-n=\frac{\sqrt{n^2+4}-n}2=[0,\bar n]\in (0,1)$ is  an isolated point of $D_{\g,\tau}\cap[0,1]$.

\Giu
As one may expect, proofs make use of the
theory of continued fractions; see \cite{HW} for general information. Let $\g>0$, $\t\ge 1$; let $\x$ be an irrational number and let: 
$$\x=[a_0,a_1,a_2,...]=a_0+\su{a_1+\su{a_2+\su{\ddots}}}$$ 
be its continued fraction expansion, $p_k/q_k=[a_0,a_1,...,a_k]$ its $k^{\rm th}$ convergent, and $a'_k:=[a_k,a_{k+1},...]$ its $k^{\rm th}$ complete quotient.  

\Giu
{\bf Lemma}
{\sl
A number $\x$ belongs to $D_{\g,\t}$ if and only if
\beq{condition}
\frac{q_{k+1}}{q_{k}^{\t}}+\frac{1}{a'_{k+2}q_{k}^{\t-1}}\leq\frac{1}{\g}\ ,\qquad \forall k\ge 0\ .\eeq
}

\Giu
{\bf Proof}
From continued fraction theory, one knows  that\footnote{For a proof, see, e.g., Lemma 1, Appendix 8, p. 122 of \cite{CC}.}
\beq{lemmetto}
\x\in D_{\g,\t} \quad \iff \quad  \Big|\x- \frac{p_k}{q_k}\Big|\ge \frac{\g}{q_{k}^{\t+1}}\ ,\qquad \forall k\ge 0\,.
\eeq
Then\footnote{$p_{-1}:= 1$, $q_{-1}:=0$; for the first equality, see \cite[\S 10.7]{HW}.},
\beqano
\Big|\x- \frac{p_k}{q_k}\Big|&=& \frac{1}{q_{k}(a'_{k+1}q_{k}+q_{k-1})}
=\su{{q_k}^{\t+1}}\ \frac{q_k^\t}{a'_{k+1}q_k+q_{k-1}}\\
&=&\su{{q_k}^{\t+1}}\ \frac{q_k^\t}{a_{k+1}q_k+q_{k-1}+\frac{q_k}{a'_{k+2}}}\\
&=&\su{{q_k}^{\t+1}}\ \frac{q_k^\t}{q_{k+1}+\frac{q_k}{a'_{k+2}}} \\
&=& \su{{q_k}^{\t+1}}\  \Big( \frac{q_{k+1}}{q_k^\t}+\frac{1}{a'_{k+2} q_k^{\t-1} }\Big)^{-1}\ ,
\eeqano
and the claim follows from \equ{lemmetto}. \qed

\Giu
{\bf Proof} (of the Theorem~1)  One immediately verifies that
\beq{alpha}
\left\{
\begin{array}{l}
\dst \a=n+\su{\a}\ , \qquad n^\t=\a\ ,\\  \ \\
\dst \a=[n,n,n,n,....]\ ,\\ \ \\
p_0=n,\ q_0=1,\ p_1=n^2+1,\  q_1=n,\  a'_k=\a,\ q_{k+1}=p_k\ (\forall k\ge 0)\ .
\end{array}\right.
\eeq
Thus, for $k=0$ we have
\beq{a0}
\Big| \a- \frac{p_0}{q_0}\Big|\ugper{\equ{alpha}} \a-n\ugper{\equ{alpha}} \su{\a}\ugper{\equ{gt}}\g .
\eeq
For $k\ge 1$, using \equ{alpha} and the facts that $p_k/q_k\le p_1/q_1$ and $q_k\ge q_1$, one finds
\beqano
\frac{q_{k+1}}{q_{k}^{\t}}+\frac{1}{a'_{k+2}q_{k}^{\t-1}}&\ugper{}&\frac{p_{k}}{q_{k}}\frac{1}{q_{k}^{\t-1}}+\frac{1}{\a q_{k}^{\t-1}}\\
&\leq& \frac{p_{1}}{q_{1}}\frac{1}{q_{1}^{\t-1}}+\frac{1}{\a q_{1}^{\t-1}}\ugper{}\frac{n^2+1}{n^\t}+\su{n^{\t-1}\a}\\
&\ugper{}& \frac{n^2+1}{\a} +\frac{n}{\a^2}=\su{\a} \Big(n^2+1+\frac{n}{\a}\Big)\\
&=&\su{\a} \ (\a n+1)=n+\su{\a} =\a\\
&=&\su\g\ ,
\eeqano
which, together with \equ{a0} and the Lemma, shows that $\a\in D_{\g,\t}$. 
\\
Next,  because of \equ{alpha}, 
\beqano
\Big| \a -\frac{p_1}{q_1}\Big|&=& \frac{p_1}{q_1} - \a= \frac{n^2+1}{n}-\a=\su{n}+n-\a\\
&=&\su{n}-\su{\a}=\su{n\a^2}=
\su\a\ \su{q_1 n^\t}=\su{\a q_1^{\t+1}}\\
&=&\frac{\g}{q_1^{\t+1}}\,.
\eeqano
Such relation, together with \equ{a0}, shows that 
$\a$ separates the two intervals\footnote{Recall the definition of the open intervals $I_{\g,\t}(p,q)$ in \equ{compl}.}  $I_{\g,\t}(p_0,q_0)$ and $I_{\g,\t}(p_1,q_1)$, and, therefore, $\a$ is an {\sl isolated point} of $D_{\g,\t}$. \qed

\section{General  properties of isolated Diophantine numbers}

General properties of isolated Diophantine numbers have been investigated in \cite{F1,F2}, where  proofs may be found. Let us briefly report here the main results in \cite{F1,F2}.

 \nl
 The first result in \cite{F1} shows that isolated points are not that rare: Indeed, {\sl any} Diophantine number has at least one equivalent representative\footnote{Recall that two irrational numbers  $\x$ and $\x'$ are {\sl equivalent} if and only if $\x'=\frac{a \x+ b}{c \x +d}$ with integers $a,b,c,d$ satisfying $ad-bc=\pm 1$ , and that happens if and only if the continued fractions of $\x$ and $\x'$ differ only by a finite number of terms; compare \cite[\S 10.11]{HW}.}, 
 which is isolated in some Diophantine set:

\Giu
{\bf Theorem 2} \cite[Theorem B]{F1} {\sl Fix $\g\in(0,\su{2})$, $\t\geq 1$,  $\a\in D_{\g,\t}$, and   let $m:=\big[\frac{3\cdot 2^{\t} }{\g}\big]$. Then, the equivalent Diophantine number 
 $\dst \a':=\frac{m \a+1}{(2m+1)\a+2}$  is an isolated point of 
$D_{\gamma_\alpha, \tau_\alpha}$ for suitable $\t_{\a}>\t$ and $\gamma_\a>0$.}

\Giu
Actually, it can happen that a  Diophantine number is {\sl simultaneously isolated for infinitely many Diophantine sets}. More precisely\footnote{In \cite{F1} something  stronger is proven, in the sense that the sequence $\t_n$ in Theorem~3 can be assigned arbitrarily up to small errors. }:

\giu
{\bf Theorem 3} \cite[Theorem A]{F1} {\sl For all $\t\geq 1$, there exist $\g>0$ and $\a\in D_{\g,\t}$ such that $\a$ is an isolated point for $D_{\g_n,\t_{n}}$, for suitable sequences $\t_n\searrow\t$, $\g_n\searrow \g$.}

\Giu
Even though these two theorems show the existence of many isolated Diophantine numbers, 
from the metric point of view, the typical situation seems to be that Diophantine sets are Cantor sets:

\Giu
{\bf{Theorem 4 \cite{F2}}}
{\sl Let $\t>\t_0:=\frac{3+\sqrt{17}}2=[\overline{3,1,1}]$. Then, for almost all $\g\in (0,1/2)$, $D_{\g,\t}$ is a Cantor set.}

\section{Remarks}

(i)
The Diophantine exponent $\t_0:=\frac{3+\sqrt{17}}{2}$  in Theorem 4 is certainly not optimal, and it would  not be difficult to improve it. On the other hand, it is not so obvious  
 what is the {\sl optimal} $\t_0$, for which  the statement of Theorem 4 holds.

 \giu
 (ii) Diophantine sets, as pointed out in the Introduction, play a fundamental r\^ole in Dynamics, e.g., in the theory of exact symplectic twist diffeomorphisms. Arithmetic properties of the rotation number of an invariant curve of a twist diffeomorphism  are, in particular, relevant for the renormalization point of view; compare \cite{MS}. Now, even though Theorem~2 above indicates that the property of being isolated for Diophantine numbers may not be a stable property under renormalization, it would be quite interesting to see if 
 such a property  does have a counter part in dynamics. For example, 
 
 \nl
 {\sl
 Does there exists a $C^r$ exact symplectic twist diffeomorphism $f$, $r\ge 2$, having an isolated 
 invariant curve of rotation number $\a$ that is not of bounded type, with $\a$ isolated point of a suitable Diophantine set?
}
 
\giu
(iii) A final comment on higher dimensional Diophantine sets.\\
Let $n\geq1$, $\g,\t>0$, an define 
$$D^{n}_{\g,\t}:=\{\a\in\Bbb{R}^n: |q\cdot\a-p|\geq \frac{\g}{|q|^{\t}}\,, \forall q\in\Bbb{Z}^n\bks\{0\}, p\in\Bbb{Z}\}\,.
$$
The analogous problem discussed in this note is\footnote{Homogeneous Diophantine sets $\real^n_{\g,\t}:=\{\o\in \real^n: |\o\cdot k|\ge \g/|k|^\t\,, \ \forall k\in \integer^n\,, k\neq 0\}$ are trivially Cantor sets, since if $\o\in  \real^n_{\g,\t}$, then also $t\o\in \real^n_{\g,\t}$, $\forall t\ge 1$; see  \cite{Br} for related questions on  $\real^n_{\g,\t}$.}:

\nl
 {\sl For $n\ge 2$, do there exist $\g,\t>0$ such that $D^{n}_{\g,\t}$ is not a Cantor set?}  

\nl
Clearly, such a question may be more difficult to analyse due to the lack of the beautiful and powerful   theory  
of continued fractions.

\Giu
{\footnotesize 
{\bf Acknowledgements} We are indebted with Michel Waldschmidt for useful discussions and for his encouragement.  We also thank, Yann Bugeaud for interesting comments.
}

\end{document}